\documentclass[dvips,11pt]{article}

\usepackage{epsfig,}
\usepackage{amssymb}
\usepackage{fullpage}

\newtheorem{theorem}{Theorem}[section]
\newtheorem{lemma}[theorem]{Lemma}
\newtheorem{corollary}[theorem]{Corollary}
\newtheorem{proposition}[theorem]{Proposition}

\newtheorem{definition}[theorem]{Definition}

\newcommand{\Top}{{\bf t}}
\newcommand{\bottom}{{\bf b}}
\newcommand{\fcF}{{\hat{G}_F}}
\newcommand{\fcvF}{{\hat{V}_F}}
\newcommand{\fceF}{{\hat{E}_F}}

\makeatletter
\def\Box{\raise3pt\hbox{\fbox{\kern0.07cm}}}
\def\@yproof[#1]{\@proof{ #1}}
\def\@proof#1{\begin{trivlist}\item[]{\em Proof#1.}}
\newenvironment{proof}{\@ifnextchar[{\@yproof}{\@proof{} 
}}{\hfill$\Box$\end{trivlist}\makeatother}

\newcommand{\N}{\mathbb{N}}
\newcommand{\Z}{\mathbb{Z}}
\newcommand{\R}{\mathbb{R}}

\begin{document}

\thispagestyle{empty}    
    
\title{Domino tilings  and  related models: \\
 space of configurations of domains with holes}

\author{S\'ebastien Desreux
\thanks{Laboratoire d'Informatique Algorithmique: Fondements
et Applications, Universit\'e  Paris 7,
2 place Jussieu 75251 Paris cedex 05, France 
{\tt  Sebastien.Desreux@liafa.jussieu.fr}.}
\and
Martin Matamala
\thanks{Departamento de Ingenieria Matematica and Centro de Modelamiento
Matematico UMR 2071 U. Chile-CNRS, 
Blanco Encalada 2120, Santiago, Chile
{\tt mmatamal@dim.uchile.cl}.}
\and
Ivan Rapaport
\thanks{Departamento de Ingenieria Matematica and Centro de Modelamiento
Matematico UMR 2071 U. Chile-CNRS,
Blanco Encalada 2120, Santiago, Chile
{\tt irapapor@dim.uchile.cl}.}
\and 
Eric R\'emila
\thanks{Laboratoire de l'Informatique du Parall\'elisme,
Unit\'e Mixte de Recherche  5668, ENS Lyon-CNRS-INRIA,
46 All\'ee d'Italie, 69364 Lyon Cedex 07, France 
{\em and}
Groupe de Recherche en Informatique et Math\'ematiques Appliqu\'ees,
IUT Roanne (Univ. St-Etienne), 20 avenue de Paris, 
42334 Roanne Cedex, France 
{\tt Eric.Remila@ens-lyon.fr}.}
}

\maketitle

\begin{abstract}
 
We first prove that the set of domino tilings of a fixed finite figure is a 
distributive lattice, even in the case when the figure has holes. 
We then give a geometrical interpretation of the order given by 
this lattice, using (not necessarily local) transformations called {\em flips}.

This study allows us to formulate an exhaustive generation algorithm and a
uniform random sampling algorithm. 

We finally extend these results to other types of tilings (calisson tilings,
tilings with bicolored Wang tiles). 

\end{abstract}

\section{Introduction} 

In the last ten years, a lot of progress has been done about 
the study of tilings. Most remarkably, W. P. Thurston \cite{Thu}, using work 
of J. H. Conway and J. F. Lagarias \cite{Con-Lag}, introduced the 
notion of height functions,  which encode domino tilings and 
calisson tilings of a polygon $P$. 

The notion of height function appears to be a very powerful tool for the
study of tilings. It has notably been extended by different authors 
\cite{Ken-Ken} \cite{Rem1} to study  tiling algorithms for other 
sets of prototiles.

For domino  tilings, height functions induce a lattice 
structure on the set of tilings of a fixed polygon (see \cite{Rem_Ordal}). 
Some important results are obtained from this 
structure: A linear time 
tiling algorithm \cite{Thu}, rapidly mixing Markov chains for 
random sampling 
\cite{lubyrandallsinclair} \cite{Wil}, computation of the number of 
necessary flips (local 
transformations involving two dominoes) to pass 
from a fixed tiling to another fixed tiling \cite{Rem_Ordal}, 
efficient exhaustive generation of tilings \cite{Desreux} 
\cite{Des-Rem}.   

Dominoes  are of particular importance 
to theoretical physicists, for 
whom dominoes  are models of \emph{dimers}, which are diatomic molecules 
(such as dihydrogen), and each tiling is seen as a possible state of 
a solid or a fluid. 

The present paper tries to generalize previous results  to  figures 
which are not polygons, {\it i. e.}  figures with holes. This is done by the 
introduction of an {\em equilibrium function} on edges of cells of 
the figure. With this tool, we  prove that the set of tilings of any 
finite figure has a 
distributive lattice structure, of which we give a geometrical 
interpretation
of this structure. To this end, 
we also need to introduce some structural notions, following works of J. 
Propp \cite{Prop} and J. C. Fournier \cite{Fournier}: The   
{\em critical cycles}, which induce {\em forced components} and {\em  
generalized flips}.

Our approach is constructive; this allows us to 
exhibit algorithms  to 
compute the objects introduced. As a consequence, we obtain an exhaustive 
generation algorithm and a uniform random sampling algorithm.

We finish by proving that these ideas can  be directly adapted for other 
types of tilings: Calisson tilings and tilings with bicolored Wang 
tiles.

\section{Figures in the plane grid}

\subsection{The plane grid} 

Let $\Lambda$ be the plane grid of the Euclidean plane $\R^2$. A 
{\em vertex} of $\Lambda$ is a point whose coordinates are both integers.

A vertex $v=(x_1,y_1)$ is a  \emph{neighbour} of another vertex $v=(x_2,y_2)$ 
if $|x_1-x_2|+|y_1-y_2|=1$. Hence, each vertex $v$ has four 
neighbours $v+(1,0)$, $v-(1,0)$, $v+(0,1)$ and
$v-(0,1)$ 
which are canonically called the {\em  East, West, North} and 
{\em South neighbour} of $v$, respectively. 
An \emph{edge} of $\Lambda$ 
is the closed segment of straight line between two adjacent vertices.
A {\em cell} of $\Lambda$ is a (closed) unit square whose 
corners are  vertices. 
Two cells are {\em 4-neighbours} (respectively {\em 8-neighbours}) 
if they share an edge (respectively at least a vertex).

A directed graph $G=(V,E)$ is \emph{symmetric} if $(v,v')\in E$
if and only if $(v',v)\in E$ for all $v,v'\in V$.
In this work we deal with the symmetric directed graph 
(denoted by $\Lambda^+$) obtained from the planar
grid $\Lambda$ by replacing  each edge $vv'$ by
two \emph{arcs} $(v,v')$ and $(v',v)$. For an 
arc $a=(v,v')$ of $\Lambda^+$ we denote by $[a]$
its associated edge in $\Lambda$.

A \emph{(directed) path} $P$ in a directed graph $G=(V,E)$ 
is a sequence of vertices $(v_0,\ldots,v_k)$
such that $(v_i,v_{i+1})$ is an arc of $G$
for every $i=0,\ldots,k-1$.
We denote by $E(P)$ the multiset of
all the arcs used by the path $P$ and by $V(P)$ the multiset of its vertices.
We say that $G$ is \emph{connected} if any two vertices  
of $V$ are linked by a path.

A  path $P = (v_0,\ldots,v_k)$ with $v_k = v_0$ is called 
a {\em cycle}.
The cycle is {\em elementary} if $v_{i} = v_{j}$ and $i \neq 
j$ imply $ \{ i, j \} =  \{ 0, k \}$. 
In a plane graph one has two kinds of elementary  cycles:
The \emph{clockwise} cycles and the \emph{counterclockwise} ones.  


Let $G=(V,E)$ be a symmetric directed graph. A function 
$g:E\to \mathbb{Z}$ is \emph{skew-symmetric}
if $g(v,v')=-g(v',v)$ for every $(v,v')\in E$.
Given any function $h:V\to \mathbb{Z}$
we  define its associated difference function 
$D(h):E\to \mathbb{Z}$ by $D(h)(v,v')=h(v')-h(v)$, 
for all $(v,v')\in E$.
Conversely, if $G$ is connected, given a function $g:E\to \mathbb{Z}$
which satisfies $g(C)=0$ for all cycle $C$ of $G_F$ 
and a vertex $w_0$ of $V$, there exists a unique function
$h:V\to \mathbb{Z}$ such that $h(w_0)=0$ and 
$D(h)(a)=g(a)$ for all $a\in E$.

Let $E'$ be a multiset of arcs of $G$. 
We denote by $g(E')$ the sum of the values $g(a)$
over all the arcs $a\in E'$ (each arc $a$ is counted according to its 
corresponding multiplicity in the multiset).
Then $g(E')=\sum_{a\in E'}g(a)$. 
For a path $P$, instead of $g(E(P))$ we use the shorthand $g(P)$.

We assume that cells of $\Lambda$ are colored as a checkerboard. 
We thus  have black cells and white cells, and two cells sharing 
an edge have different colors. 
Let us define the \emph{spin} function $sp$ on the arcs of $\Lambda^+$. 
For each arc $a=(v, v')$, the spin of $a$ is
noted $sp(a)$ and given by:
\begin{itemize} 
\item $sp(a) = 1$  if an ant moving from $v$ to $v'$  
has a white cell on its left side (and a black cell on its right 
side);
\item $sp(a) = -1$ otherwise. 
\end{itemize} 

For each clockwise elementary cycle 
$C$, one has  $sp(C)=4Dis(C)$ where 
$Dis(C)$, the \emph{disequilibrium} of $C$, 
is the difference between the number of black cells 
and the number of white cells enclosed by $C$. 
The result is true for each cycle enclosing a single cell, 
and each elementary cycle can be decomposed into such sqaure cycles.

\subsection{Figures} 

A {\em figure} $F$ of $\Lambda$ is a 4-connected, finite union of cells
of $\Lambda$.
The unique infinite 8-connected 
component of $\R^2 \setminus F$ is denoted by $H_{\infty}$. 
 The other ones 
are called the {\em holes} of $F$. 
The set of all edges in the boundary 
of $F$ is denoted by $E_b(F)$. 
The set of edges in $H_{\infty}\cap F$ is called the 
\emph{outer-boundary} of $F$ and denoted by $E_{ob}(F)$.
Analogously, we denote by $V_b(F)$ the set 
of all the vertices  on the boundary of $F$.

\begin{figure}[htbp]
\centerline{\epsfig{file=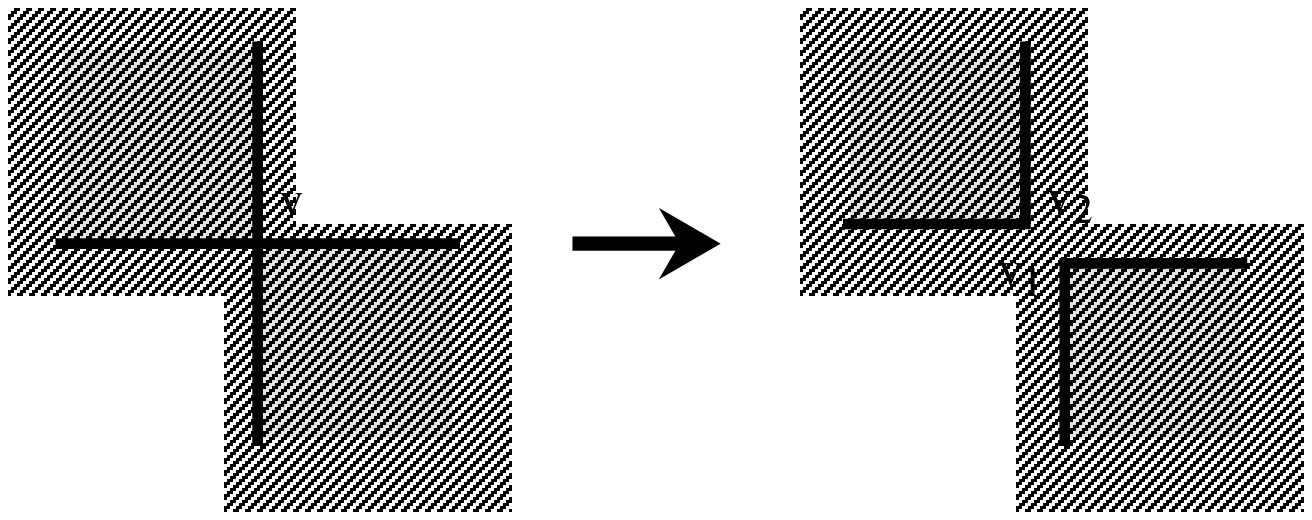,width=6cm}}
\caption{Vertex duplication according to 4-connectivity of $F$ 
and 8-connectivity of $\R^2 \setminus F$.  
} 
\label{vertex_duplic}
\end{figure}

Because of the two types of
connectivity 
for cells, we 
replace (until the end of the paper)  each vertex $v$ of $F$ 
such that each edge issued from $v$ is 
on the boundary of $F$, by two vertices $v_{1}$ and $v_{2}$, each of 
them connected to exactly two neighbours of $v$ (see Figure
\ref{vertex_duplic}).

A figure defines a symmetric directed 
graph $G_{F} = (V_{F}, E_{F} )$ such that 
$V_{F}$ is the set of corners of cells of $F$ (once duplication is 
done), and  $E_{F}$ is the set of arcs $a$ 
such that $[a]$ is a side of a cell of $F$. 
From this point of view, the clockwise and counterclockwise contours
of each hole are elementary cycles of $G_{F}$.

For each elementary clockwise cycle of $F$ (\emph{i. e.}  whose arcs are in 
$E_{F}$) we define  $Dis_{F}(C)$ as the difference between 
the number of black cells of $F$ 
and the number of white cells of $F$ enclosed by $C$.

\subsection{Equilibrium function} 

Informally, we can say that we want to work as if $F$ had no hole. 
To this end, the informal idea is to introduce values on edges which make 
holes disappear. 
Precisely,  this is done by the use of {\em equilibrium functions} as defined 
below: 

\begin{definition}  
An equilibrium function (denoted by $eq$) is a 
skew-symmetric function from $E_{F}$
    to $\Z$ such that
$sp(C) +eq(C)=4Dis_F(C)$, for every clockwise cycle $C$ of $F$.
\end {definition}

For figures without holes, it suffices to take  $eq = 0$.

Notice that, through the decomposition of cycles,  
a function $eq$ is an equilibrium
function if and only if   the following conditions holds.  
\begin{itemize} 
\item $eq(C)= 0$ for each elementary cycle $C$ around a cell of $F$.
\item $eq(C)=-sp(C)$ for each cycle $C$ 
which follows clockwise the boundary of a hole of $F$.
\end{itemize} 
More generally, for every  cycle $C$ of $F$, $eq(C)$ is a number 
which does not depend on the chosen equilibrium value. 
We prove in Section \ref{s:cons-alg} that 
every figure has an equilibrium function
which can be efficiently computed.

We also need some auxiliary functions deduced from the function $eq$.

\begin{definition}
The functions $eq_r$, $\Top$ and $\bottom$ are defined as follows.
\begin{enumerate}
\item $eq_r(a)=eq(a)-sp(a)$ for all $a\in E_F$.
\item $\Top(a)=eq(a)-sp(a)+2$ for all $a\in E_F\setminus E_b(F)$.
\item $\bottom(a)=eq(a)-sp(a)-2$ for all $a\in E_F\setminus E_b(F)$.
\item $\Top(a)=\bottom(a)=eq(a)+sp(a)$ for all $a\in E_b(F)$.
\end{enumerate}
\end{definition}

Note that for any arc $a$ in $ E_F$, $\Top(a)-\bottom(a)$ 
is either 0 or 4.


\section{The lattice of tilings}

In this section we associate three classes of objects to a figure:
Tilings, height functions and acyclic orientations. Our goal is the 
study of tilings, and height functions and acyclic orientations are 
some powerful tools to do this study.

\subsection{Tilings and height functions}

A {\em domino} is a figure formed by two cells sharing an edge, 
which is called the {\em central axis} of the domino. 
A {\em tiling} $T$ of a figure $F$ is a set of dominoes included in 
$F$, with pairwise disjoint interiors (\emph{i. e.} there is no overlap), 
such that the union 
of the tiles of $T$ equals $F$ (\emph{i. e.} there is  no gap).
Each tiling $T$ of a figure $F$ is completely determined
by the set of its central axis. 
The characteristic function of a tiling $T$ defined 
from the set of arcs of $F$ into $\mathbb{Z}$ 
is given by $\chi_T(a)=1$ if $[a]$ is  a central axis of a 
tile of $T$ and $\chi_T(a)=0$ otherwise.

Let $T$ be a tiling of $F$. The {\em  height difference} in $T$,
noted $g_T$\,, is the skew-symmetric function defined by

$$\forall a\in E_F \quad g_T(a)=eq_r(a)+2sp(a)(1-2\chi_T(a))$$
 
Then $g_T(a)\in \{\bottom(a),\Top(a)\}$ for every $a\in E_F$.

Let us define $g_F(a)$ as $eq(a)+sp(a)$ for every $a\in E_F$. 
It can be seen that for each pair $(T, T')$ of tilings of $F$ and
every $a\in E_F$,
$g_T(a)-g_{T'}(a)=4sp(a)(\chi_T(a)-\chi_{T'}(a))$.
Thus, if $g_T = g_{T'}$,  then  $T = T'$. Hence,
the function $g_T$ is a tool to encode the tilings.
Moreover, for every  arc $a\in E_b(F)$ and every tiling $T$ of $F$  
we necessarily have $g_T(a)=g_F(a)$. Thus 
$g_T(a)$ does not depend on $T$ for $a\in E_b(F)$. 
Additionally, $g_T(a)-g_{T'}(a)\in \{-4,0,4\}$, 
for every $a\in E_F$.

\begin{proposition} 
    \label{coherence des hauteurs}
Let $T$ be a  tiling of a figure $F$. For each cycle $C$ of $F$, one
has $g_{T}(C) = 0$. 
\end{proposition}

This proposition is a generalization of a theorem by J. H. Conway 
\cite{Con-Lag} about tilings of polygons. 

\begin{proof} 
(sketch) It suffices to prove the result for elementary cycles since the 
height difference of each cycle is the sum of the height differences 
of the elementary cycles which compose it. This is done by induction 
on the number of cells of $\Lambda$ enclosed by the cycle.

The case of a cycle following the boundary of a 
hole is 
easily treated from the definition of equilibrium functions. 
We also verify that  the proposition holds for elementary cycles of 
length 4 around a cell. 

We now use induction. If we are not in one of the cases 
treated above, then the area enclosed by the cycle can be cut by a 
path in $F$, which induces two new cycles, each of them enclosing 
less cells of $\Lambda$ than the original cycle. Thus, by the induction 
hypothesis, the height difference of both induced cycles is null, 
from which it is easily deduced that the height difference of the 
original cycle is null. 
\end{proof}

Proposition 
\ref{coherence des hauteurs} 
guarantees the correctness of the definition below. 

\begin{definition}
For each tiling $T$, the height function induced by $T$ (denoted by 
$h_T$) is the function from the set $V_{F}$ of vertices of cells of $F$ 
(once necessary vertex duplications have been done) to 
the set $\Z$ of integers, defined by  $h_T(w_0) = 0$ 
and $D(h_{T}) = g_{T}$. 
\end{definition}

We now give a characterization of height functions of tilings. 

\begin{proposition}\label{caracterisation_des_hauteurs}

Let $w_0$ be a fixed vertex of $H_\infty$. 
We denote by $\mathcal{H}_F$ the class of all the 
 functions 
 $h:V_F\to \Z$ satisfying the following properties:
\begin{itemize}
\item $h(w_0) =  0$;
\item for each arc $a$ of $E_{F}$, $D(h)(a)\in\{\bottom(a),\Top(a)\}$.
\end{itemize}

For each tiling $T$, the function $h_T$ belongs to $\mathcal{H}_F$.
Conversely, for each $h\in \mathcal{H}_F$ there exists
a tiling $T$ such that $h=h_T$. 
\end{proposition}

\begin{proof}
The first statement follows directly from the definition of $g_T$. 
Let $h$ belong to $\mathcal{H}_T$ and  let $C$ be a cycle around a cell. 
Then $D(h)(C)-eq_r(C)=sp(C)$ since $D(h)(C)=eq(C)=0$.  
Clearly, $|sp(C)|=4$. Since $D(h)(a)\in \{\bottom(a),\Top(a)\}$,
there are three arcs of $C$ such that $D(h)(a)-eq_r(a)=2sp(a)$
and exactly one arc satisfying $D(h)(a)-eq_r(a)=-2sp(a)$. 
Thus, the set $T$ of all dominoes whose 
central axis $[a]$ is such that $hd(a)-eq_r(a)=-2sp(a)$
is a tiling $T$ of $F$. 
The equality $h(v) = h_T(v)$, for each vertex $v$ of $F$, is obvious
by induction on the distance from $w_0$ to $v$.
\end{proof}

The proposition above allows one to consider each tiling as a height 
function. 

\begin{lemma}\label{cong}
For every pair of height functions $h$ and $h'$ 
and for each vertex $v$ of $F$, one has 
$h(v) - h'(v) = 0[4]$. Moreover, $h(v) - h'(v)$ does not 
depend on the chosen equilibrium function. 
\end{lemma}
\begin{proof}
Obvious by induction on the length of a shortest path from $w_0$ to 
$v$. 
\end{proof}

Let $\leq$ be the canonical order on functions:
$h\leq h'$ if and only if $h(v)\leq h(v')$ for all $v\in V_F$.

\begin{proposition}\label{p:min-max-height-inv}
Let $h$ and $h'$ belong to $\mathcal{H}_F$.
The functions $inf(h,h')$ and $sup(h,h')$ belong to $\mathcal{H}_F$.
\end{proposition}

In the vocabulary of order theory (see for example \cite{Birk}, \cite{Davey}) 
the above proposition can be restated as follows: 
 $(\mathcal{H}_F,\leq)$ is a distributive lattice.

\begin{proof}
Let $h_1=inf(h,h')$. We shall prove that for every arc $(v,v')\in E_F$,
$h_1(v')-h_1(v)\in \{h(v')-h(v),h'(v')-h'(v)\}$ (the proof for $sup(h,h')$
is similar).

For the sake of contradiction, let us assume that 
there exists an arc $a=(v,v')$ of $F$ such
that $h_1(v)=h(v)<h'(v)$ and $h_1(v')=h'(v')<h(v')$.
From Lemma \ref{cong}, one has
$h'(v') \leq h(v')-4$ and $h(v) \leq h'(v)-4$.
Then $\alpha:=h'(v')-h(v')+h(v)-h'(v)$ satisfies 
$\alpha\leq -8$.
On the other hand, $\alpha=h'(v')-h'(v)-(h(v')-h(v))$.
Since $h$ and $h'$ belong to $\mathcal{H}_F$, one obtains
$\alpha \geq \bottom(v,v')-\Top(v,v')=-4$
which contradicts the hypothesis.
\end{proof}


We define the following order on $\Gamma_F$: $T\leq T'$ if
$h_T \leq h_{T'}$.
From Proposition \ref{caracterisation_des_hauteurs}, $(\Gamma_F,\leq)$ 
is isomorphic to $(\mathcal{H}_F,\leq)$.
Thus 
$(\Gamma_F,\leq)$ is a distributive lattice.

\subsection{Forced components}

\begin{definition}
   An  elementary cycle   $C$ of $G_F$  
is {\em  critical}  if  $\Top(C) = 0$. 
 
Such a cycle is {\em  strongly critical} if, moreover, 
$sp(a)=1$, for each arc $a\in E(C)\setminus E_b(F)$.

We say that $v$ and $v'$ are \emph{critically equivalent}
if there is exists 
a critical cycle $C$ such that $v,v'\in V(C)$.
The equivalence classes of this equivalence relation
are called the \emph{forced components} of the figure $F$.

Let $\fcF = (\fcvF, \fceF)$ be the symmetric graph 
whose vertices are the forced components and 
where  $(U,U')$ is an arc of  $\fcF$  if 
there exists $v\in U$ and $v'\in U'$ such that
$(v,v')$ is an arc of $G_F$. This graph is called the \emph{graph of 
forced components} of $F$. 
\end{definition} 

Notice that each boundary cycle of a hole of $F$ is 
strongly critical. The boundary cycle of the outer boundary is critical 
if and only if $F$ contains as many black cells as white cells. 
A strongly critical cycle can be deduced from each critical cycle 
by replacing each interior arc $a$ with $sp(a)=-1$ 
by a sequence $a',a'',a'''$
of three arcs of positive spins.

\begin{figure}[htbp]
\centerline{\epsfig{file=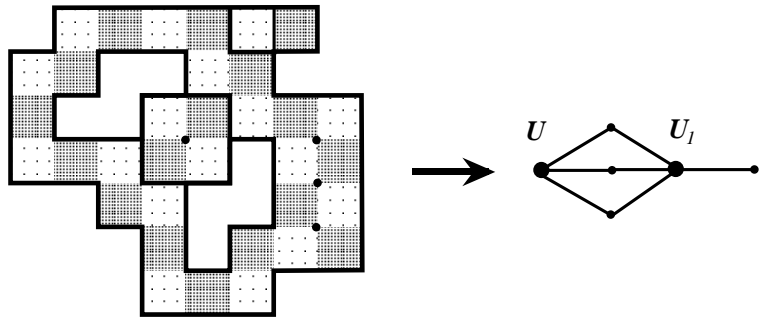,width=10cm}}
\caption{The forced components of a figure and the graph of forced 
components. 
($U_{1}$ denotes the component containing the contours of the two
holes)}
\label{forced_component}
\end{figure}

It is easy to see that for any  cycle $C$ of $G_F$, $\Top(C)$ and
$eq(C)$ do not 
depend on the chosen equilibrium  function. Thus, the notion of critical
 cycle only depends on the shape of the figure.  

\begin{definition}
Let $T$ be a tiling of $F$. 
The graph of $T$ is the 
spanning subgraph of $G_F$, denoted by $G_T=(V_F,E')$, 
where $a\in E'$ if and only if $g_T(a)=\Top(a)$.
\end{definition}  

By definition, every arc in $E_b(F)$ is an arc of $G_T$.
Moreover, $a\in E_T\setminus E_b(F)$ if and only if 
$\chi_T(a)=0$ and $sp(a)=1$, or $\chi_T(a)=1$ and
$sp(a)=-1$.

The following proposition is the reason why we are interested in critical 
cycles and forced 
components.

\begin{proposition}
    \label{cycle critique}
Let $C$ be a cycle of $G_F$. 
\begin{itemize}
\item If $C$ is critical, then for every tiling $T$ 
the cycle $C$ is a cycle of $G_T$. Conversely, if
$C$ is a cycle of $G_T$ for some tiling $T$, then 
$C$ is critical.
\item If $C$ is strongly critical, then for every tiling
$T$ the cycle $C$ is a cycle of $G_{T}$ which does not cut any tile of $T$. 
\end{itemize}
\end{proposition}

\begin{proof}
If $C$ satisfies $\Top(C)=0$, then for every tiling $T$ of $F$ 
one has $\Top(C)=g_T(C)$. Therefore
$\Top(a)=g_T(a)$ for every $a\in E(C)$, whence $C$ 
is a cycle of $G_T$.
Conversely, if $C$ is a cycle of $G_T$, then 
$g_T(a)=\Top(a)$ for every $a\in E(C)$. Thus, 
from Proposition \ref{coherence des hauteurs}, $\Top(C)=g_T(C)=0$ 
and $C$ is a critical cycle.

For the second part, if $C$ is strongly critical 
then  $\Top(a)>g_T(a)$ implies $sp(a)=1$, from which one knows
that $\chi_T(a)=0$; this means that $C$ does not cut any tile of $T$.
Conversely, let $[a]$ be an interior edge which does not cut
any tile of $T$. By definition, 
$g_{T}(a) = sp(a)+eq(a)$. Moreover,
since, $a$ belongs to $E_T$, one has
$g_{T}(a) = \Top(a)=eq(a)-sp(a)+2$ and finally $sp(a) = 1$. 
\end{proof}

\begin{corollary}
    \label{coro cycle critique}
    If $F$ has a strongly critical cycle $  (v_0, v_1, \ldots, 
v_p)$ such that for each integer $0 \leq i < p$,
$[v_i, v_{i+1}]$ an interior edge, then there exists no tiling of 
$F$. 
\end{corollary}

\begin{proof}
    Let $v_j = (x_j, y_j)$ be the vertex of this 
cycle with $x_j+y_j$ maximal, and, moreover, $x_j$ minimal with respect to the 
previous condition. One necessarily has $v_{j-1} = v_j + (-1, 0)$, 
$v_{j-2} = v_{j-1}+(0, -1)$ and (moreover) $v_{j+1} = v_j+(0,-1)$. 
Now, follow the cycle until a vertex $v_{j+2k}$ such that $v_{j+2k} 
\neq v_j+ (k,-k)$ (see Figure \ref{chemin_critique}).

\begin{figure}[htbp]
\centerline{\epsfig{file=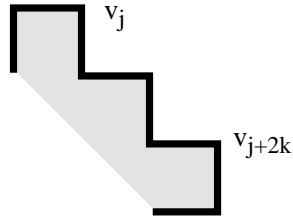,height=3cm}}
\caption{Proof of Corollary  \ref{coro cycle critique}}
\label{chemin_critique}
\end{figure}

Let $T$ be a tiling of $F$; at least one tile of $T$ must be cut by 
an edge of the path formed from the part of the cycle from  $v_j$ to 
$v_{j+2k}$.   But this is impossible, from  Proposition \ref{cycle critique}. 
Thus there exists no tiling. 
\end{proof}

From Corollary  \ref{coro cycle critique}, if $F$ can be tiled then there
are three kinds of 
forced components: The component    $U_{\infty}$, 
 the {\em single components} which are reduced to a single vertex, and the  
 {\em hole components} which contain the contour of at least one hole.

The following lemma establish an useful relation
between the height functions and the forced components.

\begin{lemma}\label{distance-equivalence}
Let $v$ and $v'$ be critically equivalent vertices. 
For all $h$ and $h'$ in $\mathcal{H}_F$,
one has $h(v) - h'(v) = h(v') - h'(v')$. 
\end{lemma}

\begin{proof}
    Let $C = (v_{0}, v_{1}, \ldots, v_{p})$ be a critical cycle
    passing through $v$ and $v'$. One can assume without loss of
    generality that $v = v_{0}$ and $v' = v_{k}$. One has
    $h(v') -h(v) = \sum_{i =0}^{k-1}D(h)(v_{i}),(v_{i+1}))
=\sum_{i =0}^{k-1}\Top(v_{i},v_{i+1}) = 
    h'(v') - h'(v)$, which yields the result. 
\end{proof}

Let us choose one vertex $v_U$ in each forced component
$U$ of $G_F$. 
From Lemma \ref{distance-equivalence},
$h\leq h'$ if and only if $h(v_U)\leq h'(v_U)$
for all $U\in \fcvF$.

We define a distance on $\mathcal{H}_F$ by
$$\Delta(h,h'):=\sum_{U\in \fcvF}\left|h(v_U)-h'(v_U)\right|$$

Notice that the distance satisfies the following equalities: 
$$\begin{array}{l}
\Delta(h,h') = \Delta(h,inf(h,h'))+ \Delta(inf(h, h'),h') \\
\Delta(h,h') = \Delta(h,sup(h,h'))+ \Delta(sup(h, h'),h') \\
\end{array}$$

\subsection{Acyclic orientations and flips}

In this part we prove that 
the lattice of height functions of a tileable figure $F$
is isomorphic to a lattice of a subclass of
orientations of the graph of forced components.
These lattices have been precisely studied by 
J. Propp  \cite{Prop}.

\begin{definition}
A directed graph $G$ is an  
\emph{orientation of $F$} if $G$ is an  orientation
of $\fcF$ such that 
$|E(C)\cap E(G)|=-\frac{1}{4}\bottom(C)$ for all
cycle $C$ of $\fcF$.
We denote by $\mathcal{G}_F$ the class of all 
the acyclic orientations of $F$.

Let $h$ be a height function and let $T$ be the corresponding tiling.
Let $G_h=(\fcvF,E)$ be the \emph{graph of 
strongly connected components} of $G_{T}$, \emph{i. e.} the directed 
graph defined by 
$(U,U')\in E$ if and only if there exists an arc of $G_{T}$ from a 
vertex of $U$ to a vertex of $U'$.  
\end{definition}

It is well known (see \cite{Ber} for example) 
that the graph of strongly components of any directed graph is acyclic.
Hence, 
$G_h$ is an acyclic orientation of $\fcF$.

\begin{proposition}
Let $h$ be an element of $\mathcal{H}$. Then $G_h$ belongs to $\mathcal{G}_F$.
Conversely, for each $G=(V,E)\in \mathcal{G}_F$ there
exists $h\in \mathcal{H}$ such that $G=G_h$.
\end{proposition}

\begin{proof}
Let $C$ be a cycle of $\fcF$. Then 
$0=D(h)(C)=D(h)(E(C)\cap E)+D(h)(E(C)\setminus E)$.
Since 
$D(h)(E(C)\cap E)=\Top(E(C)\cap E)$ 
and, for any arc $a$ of  $E_F$,
$\Top(a)=\bottom(a)+4$, 
one obtains: 
$D(h)(E(C)\cap E)=\bottom(E(C)\cap E)+4|E(C)\cap E|$. 
Moreover, from the definition of $G_h$ one obtains
$D(h)(E(C)\setminus E)=\bottom(E(C)\setminus E)$.
Finally,
$ 0 = 4|E(C)\cap E| +\bottom(C)$. 

Conversely, let $G=(V,E)$ belong to $\mathcal{G}_F$. 
Let $g$ be the function defined by $g(a)=\Top(a)$ if
$a\in  E$ and $g(a)=\bottom(a)$ if $a\in \fceF\setminus E$.
We prove that $g(C)=0$ for all cycle $C$. Clearly
$g(C)=g(E(C)\cap E)+g(E(C)\setminus E)$. By definition of $g$,
one has
$g(E(C)\cap E)=\Top(E(C)\cap E)$ 
and 
$g(E(C)\setminus E)=\bottom(E(C)\setminus E)$. 
Since for any arc $a$ of  $E_F$ 
$\Top(a)=\bottom(a)+4$, 
 one obtains that 
$g(E(C)\cap E)=\bottom(E(C)\cap E)+4|E(C)\cap E|$.  
Then 
$g(C)=\bottom(E(C)\cap E)+4|E(C)\cap E|+\bottom(E(C)\setminus E)
=\bottom(C)+4|E(C)\cap E|=0$. 
Thus  there exists $h\in \mathcal{H}$ such that
$h(w_0)=0$ and $D(h)=g$. 
\end{proof}



\begin{definition}
Let $G=(V,E)$ belong to $\mathcal{G}_F$ and 
let $U\neq U_\infty$ be in $V$ without 
incoming (resp. outgoing) arcs. The graph obtained from $G$ by 
an \emph{upward} (resp. \emph{downward}) \emph{flip in $U$}  
is the  acyclic directed graph
$G^{U}=(V,E^+)$ (resp. $G^{U}=(V,E^-)$ ) where 
$$E^+=E\setminus\{(U,U'): (U, U')\in E\}\cup \{(U',U):(U,U') \in E\}$$
and
$$E^-=E\setminus\{(U',U): (U', U) \in E\}\cup \{(U,U'):(U', U) \in E\}$$
\end{definition}

An  upward or downward flip in $U$ 
corresponds to reversing all the arcs incident to $U$.
If $U$ is reduced to a single vertex $v$, the flip is said {\em local}.  
If $U$ contains  the contour of a hole, we say that it is a {\em hole flip}.

\begin{proposition}
For any $G=(V,E)\in\mathcal{G}_F$, $G^{ U}=(V,E')$ belongs to $\mathcal{G}_F$.
\end{proposition}

\begin{proof}
Let $C$ be a cycle of $\fcF$. Since 
the number of arcs $(U,U')$ used by $C$ 
is equal to the number of arcs $(U',U)$ used by $C$,
one has $|E(C)\cap E| =|E(C)\cap E'|$. 
\end{proof}

Let us denote by $h^U$ the unique function in $\mathcal{H}_F$
such that $G_{h^U}=G^U$. We have $h^U = h$ in $V \setminus U$ and 
 $|h^U - h| = 4$ in $U$. The two corresponding tilings 
differ only around $U$. 
 
\begin{figure}[htbp]
\centerline{\epsfig{file=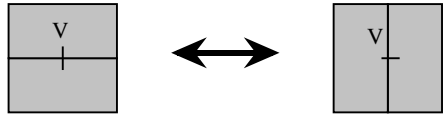,height=15mm}}
\caption{A local flip} 
\label{local flip}
\end{figure}

In particular, 
in $\Gamma_F$ a local flip in $U=\{v\}$
(see Figure \ref{local flip}) is the replacement in $T$ of the 
pair of dominoes which cover the $2 \times 2$ square centered in $v$  
by the other pair which can cover the same square.

The upward flips defined above canonically induce an order on the set 
$\mathcal{G}_F$. Given $G$ and $G'$ in $\mathcal{G}_F$,
we say that $G \leq_{flip^+} G'$ if and only if there exists a sequence 
$(U_0, \ldots, U_{p-1})$ of vertices of $\fcF$
and a sequence $(G_0, G_1, \ldots, G_p)$ of graphs of $\mathcal{G}_F$ 
such that $G_0 = G$, $G_p= G'$ and, 
for each integer $0 \leq i <  p-1$, 
$G_{i+1}$ is deduced form $G_i$ by an upward flip.

\begin{proposition}
    \label{p:equiv-order}
Let $h$ and $h'$ belong to $\mathcal{H}_F$. Then  $h\leq' h'$ if and only 
if $G_h\leq_{flip^+}G_{h'}$. 
Moreover, in this case, one can pass from $G_h$ to $G_{h'}$ by a sequence of 
$\Delta(h,h')/4$ flips.
\end{proposition}

\begin{proof}
  The direct part of the proposition is proved by induction on the 
quantity $\Delta(h, h')$. The result is obvious if $\Delta(h, h') = 
0$ ({\it i. e.} $h = h'$).

Now, let us assume that  $\Delta(h, h) \neq 0$ and $h \leq h'$. 
We will prove that there exists a forced component $U$ such that, for 
each vertex $v_{U}$ of $U$,  $h(v_{U}) < h'(v_{U})$ (which implies 
$h(v_{U}) \leq h'(v_{U}) - 4$) and an upward flip 
can be done from $G_{h}$ on $U$. 

Let $U_{0}$ be a component such that   $h(v_{U_{0}}) < h'(v_{U_{0}})$. If an
upward flip 
can be done from $G_{h}$ on $U_{0}$, then we are done. Otherwise, 
there exists
an arc $(v_{0}, v_1)$, with $v_{0}$ in $U_{0}$ and $v_{1}$ in another 
forced component $U_{1}$ such that  $D(h)(v_{0}, v_1) = \bottom(v_{0}, v_1)$.
From Lemma \ref{distance-equivalence}, we have: 
$h'(v_{U_{1}}) - h(v_{U_{1}}) = h'(v_{{1}})- h(v_{{1}})$ and
$h'(v_{U_{0}}) - h(v_{U_{0}}) = h'(v_{{0}})- h(v_{{0}})$. 
Moreover, 
$h'(v_{{1}})- h(v_{{1}}) = D(h')(v_{0}, v_1) + h'(v_{{0}})- 
h(v_{{0}}) - D(h)(v_{0}, v_1)$.
Thus $h'(v_{U_{1}}) - h(v_{U_{1}}) = h'(v_{U_{0}})
- h(v_{U_{0}}) +  D(h')(v_{0}, v_1) - \bottom(v_{0}, v_1)$, which yields
$h'(v_{U_{1}}) - h(v_{U_{1}}) \geq h'(v_{U_{0}}) - h(v_{U_{0}}) > 0$.

Either  an upward 
flip can be performed on $U_1$, or the same argument can be repeated 
from $U_1$ to obtain an arc  $(v_1, v_2)$ from $U_{1}$ to another 
forced component $U_{2}$, $D(h)(v_{1}, v_2) = \bottom(v_{1}, v_2)$.
By repeating the process, there are two possibilities: Either one obtains a 
forced  component  $U_i$ on which no
upward flip can be done, or an infinite sequence 
$(U_i)_{i \in \N}$ is obtained. But since $F$ is finite, the second 
possibility would imply that there exists a finite subsequence $(U_j, 
U_{j+1},\ldots, U_{k})$ such that $U_{j}= U_{k}$, which is a contradiction
since $G_{h}$ is acyclic.

We have now proved the existence of a forced component  $U$ such that an
upward
flip can be done from $G_{h}$ around $U$ to obtain a function $h^U$.
Notice 
that $\Delta(h^U, h') = \Delta(h, h')-4$. This proves that by induction
one can pass 
from $h$ to $h'$ with $ \Delta(h, h')/4$ flips. 

The direct part of the first part of the proposition is obvious.
\end{proof}

\begin{corollary}\label{c:isomorphic-height-flip}
The function $\varphi:\mathcal{H}_F\to \mathcal{G}_F$ defined
by $\varphi(h)=G_h$ is an order isomorphism between
$(\mathcal{H}_F,\leq)$ and  $(\mathcal{G}_F,\leq_{flip^+})$.
\end{corollary}

\begin{corollary}
Let $h$ and $h'$ belong to $\mathcal{H}_F$.
The number of successive flips needed to pass from $G_h$ to $G_{h'}$
is $\Delta(h,h')/4$. 
Moreover, the components $U$ on which flips are done in such a sequence 
are those such that
$h(v_{U}) \neq h'(v_{U})$. 
\end{corollary}

\begin{proof}
   A flip changes $\Delta(h,h')/4$ by one unit. Thus
$\Delta(h,h')/4$ is a lower bound and the bound is reached if each 
flip lets the quantity decrease.  Thus the bound can be reached only if 
    the components $U$ on which flips are done 
are precisely those such that $h(v_{U}) \neq h'(v_{U})$.

Conversely, we have seen that for $h \leq  h'$, one can pass from 
$G_h$ to $G_{h'}$ by a sequence of $\Delta(h,h')/4$ flips. For the 
general case we use  $inf(h, h')$. Passing through $G_{inf(h, h')}$, one can 
pass from $G_{h}$ to $G_{h'}$ by a sequence of 
$\Delta(h,inf(h, h'))/4 + \Delta(inf(h, h'),h')/4 = \Delta(h,h')/4$ 
successive flips. 
\end{proof}

\begin{corollary}
For each pair of tilings $(T, T')$,  it is possible to pass from 
$T$ to $T'$  by a sequence of local flips if and only if $h_{T} = 
h_{T'}$ for all the vertices on the boundary of $F$. 
Moreover, in this case the number of local flips needed is $\sum_{v 
\in E_{F}}\left|h_{T}(v)-h_{T'}(v)\right|$.
\end{corollary}

\begin{proof}
    This is a special case of the previous corollary. 
    \end{proof}

\subsubsection{Freeness and rigidity of tilings}

Informally, our tools allow one to see what is forced and
what one has to choose in order to tile a figure.  


\begin{lemma}
    \label{minimal}
Let $h$ belong to $\mathcal{H}_F$. The function $h$ is maximal if and only if
for each  forced component $U$ 
there exists a path $P$ from $U_\infty$ to $U$
in $G_h$.

The function $h$ is minimal if and only if for each 
forced component $U$ 
there exists a path $P$ from $U$ to $U_\infty$ 
in $G_h$.
\end{lemma}

\begin{proof}
If there exists $U'$ such that there is
no path from $U_\infty$ to $U'$ in $G_h$, then by
taking a longest $P'$ which finishes in $U'$
one clearly deduces the existence of a component $U\neq U_\infty$ of $G_{h}$ 
with  no incoming arcs. 
Thus an upward flip can be done in $U$ and $h$ is not maximal. 

Conversely, 
if $h\neq h_{\max}$ then it is possible to pass from 
$h$ to $h_{\max}$ by a sequence of upward flips. 
The first component $u$ on which a flip is done in such a sequence has 
no outgoing arc. 

The proof for $h_{\min}$ is similar. 
\end{proof}

\begin{proposition}
    For each component $U$ such that $U\neq U_\infty$, 
    one has $h_{\min}(v_{U}) \neq h_{\max}(v_{U})$.
\end{proposition}

\begin{proof}
    Assume that $h_{\min}(v_{U}) = h_{\max}(v_{U})$.  Take an 
    sequence of upward flips from $h_{\min}$ to  $h_{\max}$. 
    Since no flip is done on $U$, no upward flip can be done on each 
    component  on a path of $G_{h_{\min}}$ from $U$ to  $U_\infty$. 
    Thus this path is also a path of $G_{h_{\max}}$.
    By concatenating it with a  path of $G_{h_{\max}}$ from $U_\infty$ to 
    $U$, a cycle appears in $G_{h_{\max}}$, which is a contradiction. 
\end{proof}

\begin{corollary}
    A vertex  $v$  of $V_{F}$ belongs to $U_{\infty}$ if and only if 
   $h(v) = h'(v)$ for all $h,h'\in\mathcal{H}_F$. 
   
   An arc $a$ links two vertices of the same component if and only if 
   $D(h)(a) = D(h')(a)$, for all $h,h'\in\mathcal{H}_F$.
\end{corollary}

\begin{proof}
    If $v$ belongs to $U_{\infty}$, then $h(v) - h'(v) = h(w_{0})- 
      h'(w_{0}) = 0$ from Lemma \ref{distance-equivalence}. 
      Otherwise, $h_{\min}(v) \neq h_{\max}(v)$ according to the 
      previous Proposition. 
      
      If an arc $a$ links two vertices of the same component, then 
      $D(h)(a) = D(h')(a)$ from Lemma \ref{distance-equivalence}. 
      Conversely, if an arc $a$ does not link two vertices of the same 
      component, then one of its vertices is in a component $U$ such that
      $U \neq U_\infty$. 
      Take an 
    sequence of upwards flips from $h_{\min}$ to  $h_{\max}$. A flip 
    of the sequence is done on $U$, which implies  that  there exists 
    a pair $(h, h')$ of functions such that $D(h)(a) \neq D(h')(a)$.
 \end{proof}
 

\section{Effective construction and algorithms}\label{s:cons-alg}

\subsection{forced components}

When $F$ can be tiled (which is the interesting case),  the  graph 
of forced components  can be 
constructed in  polynomial  time: Given a tiling $T$, we have to construct 
the graph of strongly connected components of $G_{T}$, which can be 
done in linear time (see \cite{Tar} or \cite{Cor-Lei}).

We know from matching theory that 
there exists a $O(n^{3/2})$ algorithm \cite{Hop-Kar} to obtain such a tiling
$T$, where $n$ denotes the area of $F$. 
N.  Thiant \cite{Thian} gives an algorithm which is 
linear in the area enclosed by $F$ ({\it i. e.}  the sum of the area of $F$ 
and the areas of the holes).

\subsection{Construction of an equilibrium function}
\label{construction} 
    
\begin{figure}[htbp]
\centerline{\epsfig{file=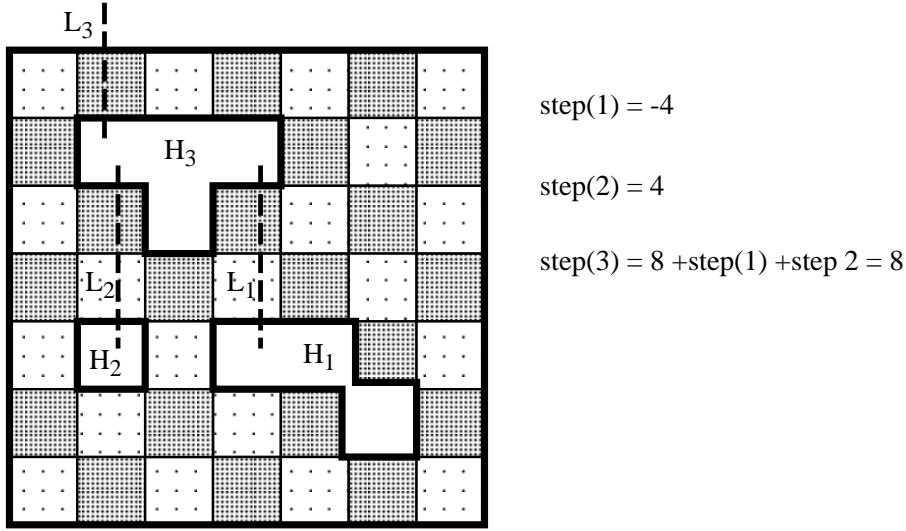,width=12cm}}
\caption{Computation of an equilibrium function by step values.} 
\label{figure}
\end{figure}

An equilibrium function can be exhibited using  {\em cut lines} (see 
also 
\cite{Sal-Tom-Cas-Rom}) as follows  
(see Figure \ref{figure}):  

For each hole $H_i$ 
of $F$, we 
(arbitrarily) fix a vertical segment $L_i = [p_i, p'_i]$ (which 
is called a {\em cut line} issued from $H_i$) of $\R^2$ such that $p_i$ is 
the central point of a highest cell of $H_i$; there exists a 
positive integer $n_i$ such that $p'_i = p_i+(0, n_i)$, the vertex 
$p'_i$ is not in $F$, and, for each integer such that $0 < n'_i < 
n_i$, the point $p_i+(0, n'_i)$ is the central point of a cell of 
$F$. Hence, the point $p'_i$ is the central point of a cell of  another 
8-connected component,
 $H_j$ of $\R^{2} \setminus F$, with $j \neq i$ (and, possibly, $j = 
\infty$).

We say that $H_j$ is the (immediate) predecessor of $H_i$. This 
construction yields a directed tree whose vertices are 8-connected 
components of $\R^2 \setminus F$. This tree is rooted in $H_{\infty}$.
We inductively define the {\em step value} of a hole 
$H_i$ (denoted by $step(i)$) as follows: Let $C_{i} = (v_{i, 0}, v_{i, 1}, 
\ldots, v_{i, p_i} )$ be a cycle which (clockwise) 
follows the boundary of $H_i$. We state: 
$$ step(i) =
\sum_{j| H_j \mbox{ has }  H_j \mbox{ for predecessor} } step(j)
-  sp(C_{i}) $$

We can now define the equilibrium value of an arc $(v, v')$ by:  

\begin{itemize} 
\item  $eq(v, v') = step (i)$ if $v'$ is the East neighbour of $v$ 
and the line segment $[v, v']$ crosses the cut line $L_i$;
\item $eq(v, v') = - step (i)$ if $v'$ is the West neighbour of $v$ 
and the line segment $[v, v']$ crosses the cut line $L_i$;
\item $eq(v, v') = 0$ if the line segment $[v, v']$ crosses no cut line. 
\end{itemize}

This function satisfies the conditions of the 
definition: It is obvious for cycles surrounding single cells, and 
if $C_{i}$ is a cycle which clockwise surrounds the hole $H_{i}$, one has
$eq(C_{i}) = step(i) - \sum_{j| H_j \mbox{ has }  H_j \mbox{ for predecessor} }
 step(j) = - sp(C_{i})$.

 The equilibrium function defined above has a specific property which can 
 be used for algorithmic arguments: For each pair $(v, v')$ of 
 vertices of $V_{F}$, there exists a path $P$ from $v$ to $v'$  such that 
 for each arc $a$ of $P$, $eq(a) = 0$. This yields that there exists 
 a spanning tree $T_{F}$ of $G_{F}$  such that for each arc $a$ of 
 $T_{F}$ (seen as a symmetric graph), $eq(a) = 0$.

 For each arc 
 $a$ of $E_{F}$, $|eq(a) |\leq 4n$ (this 
upper bound can be reached by taking a cycle which uses $a$ and cuts no 
line, except $a$). 
A precise study shows that such an equilibrium function can be 
constructed in $O(n\log n)$ time units, where $n$ denotes the number 
of cells of $F$.


\subsection{Minimal tiling}

For figures without holes, W. P. Thurston \cite{Thu} has exhibited an 
algorithm that builds the minimal tiling. For the general 
case, this algorithm can be generalized as follows.

{\bf Initialization:}
For each vertex $v$, the algorithm uses a variable value
$h(v)$ and two 
fixed values $low(v)$ and $sup(v)$. 

From previous results, we know that 
for all $a\in E_{ob}(F)$ and all tiling $T$,
$g_T(a)=g_F(a)$.
Then the value $h(v)$ can be computed for
all the vertices on the boundary of 
$H_\infty$.
If a contradiction appears, then stop. 
Otherwise, set $low(v)=sup(v):=h(v)$
for all $v$ on the boundary of $H_\infty$.

We  first construct a spanning tree $T_{F}$ on $G_{F}$
rooted in $w_{0}$. Then, for all vertices $v$ not on the 
boundary of $H_\infty$, we set 
$sup(v) = \Top(T_Fv)$ and $low(v) = \bottom(T_Fv)$,
where $T_Fv$ denotes the unique path from $w_0$ to
$v$ in $T_F$.

We set $h(v) = \bottom(v)$ for all vertices $v$ not on the 
boundary of $H_\infty$.

The algorithm also uses a set $V$ consisting of vertices $v$ such 
that there exists  a neighbour $v'$ of $v$ satisfying $h(v) + 
top(v, v') < h(v')$. 
This set is computed during the initialization.

{\bf Main loop:}
While $V$ is not empty:

\begin{itemize}
\item Pick a vertex $v$ in $V$ 
and  update $h(v)$ by adding 4 units. 
If $h(v) > sup(v)$ after updating, then stop (there is no tiling).

\item Update $V$ by adding the neighbours $v''$ of $v$ such that
$h(v'')+\Top(v'', v)
 < h(v) $ and remove $v$ if necessary. 
\end{itemize}

\begin{proposition}
    Given a figure $F$ formed of $n$ cells, 
    the above algorithm stops after at most $n^2$ time passages 
    through the loop. 
    
    Moreover, the algorithm  stops with $V$ empty if and only if there 
    exists a tiling.  In this case, when the algorithm stops, one has
    $h = h_{T_{min}}$.
\end{proposition}

\begin{proof}
    First remark that, for each arc $(v, v')$ of $E_{F}$, $h(v') - 
    h(v) = \Top(v, v')$ in 
    $\Z/4\Z$. This is 
    true during
    the initialization: We introduce the cycle $C$ formed by the concatenation 
    of (the opposite path of) $T_Fv$, $T_Fv'$ and $(v,v')$, we have
    $\Top(C) = 0[4]$ and this, together with the relations between 
    $\Top$ and $\bottom$, gives the result. 
    Moreover, this property is preserved by the loop.

    Thus each passage through the loop makes the sum 
    $\sum_{v \in V_{F}} (sup(v) - h(v))$ decrease by at least 4 
    units since, on vertex $v$, $sup(v) -low(v)
    \leq 4 |E(T_{F}v)| \leq 4n$ (where  $|E(T_{F}v)|$ denotes the number 
    of arcs of $T_{F}v$), 
    so the algorithm stops after at most 
    $n^{2}$ passages through the loop.

    When  $V$ becomes empty, each arc $(v, v')$ of $E_{F}$ 
    satisfies the hypothesis of Proposition 
    \ref{caracterisation_des_hauteurs}.  Thus there exists a tiling 
    $T$ such that  $h_{T} = h$. 
    Moreover $h \leq h_{T_{min}}$ (this is true during
    the initialization, and this property is preserved by the 
    loop), so $h_{T_{min}} = h$. 
    
    If the algorithm finds a vertex $v$ such that $h(v) > sup(v)$, then 
    there is no tiling since, otherwise, for each tiling $T$,
    $h_{T}(v) > sup(v)$, which is a contradiction (clearly, from 
    the definition of $sup$, $h_{T}(v) \leq sup(v)$). This 
    finishes the proof.  
\end{proof}

If we take into account the implementation, the algorithm's cost is at  most 
$0(n^2)$ time units, as follows:

All the values are encoded in unary numeration, which permits to add 
constant numbers in constant time.  We use the equilibrium function 
described in part \ref{construction}  
and the spanning tree $T_{F}$ is chosen in such a way that for  
each arc $a$ of $T_{F}$, $eq(a) = 0$. 
Thus,  the initialization costs $0(n^2)$ time 
units, to compute $h(v)$ for each vertex and $h(v') - h(v)$ for each 
arcs. 

 Each passage through the loop costs $O(1)$ time units since it 
    consists in a 
    fixed number of additions of 4 units and sign tests. 
    This gives the time complexity. 

Of course, a similar algorithm can be designed to construct the maximal 
tiling of $F$. 


\subsection{Exhaustive generation}

An  exhaustive generation can be done, extending ideas of \cite{Des-Rem}
to figures with holes. 
 Let $(U_{1}, U_{2}, \ldots, U_{q})$ be a fixed total order of  forced 
  components of $F$ (except $U_{\infty}$). We define a total order $<_{lex}$
on tilings 
  of $F$ as follows: Given two tilings $T$ and $T'$, we have $T 
  <_{lex} T'$ if there exists an integer $1 \leq i \leq 
  p$ such that $h_{T}(v) < h_{T'}(v)$ for each vertex $v$ of $U_{i}$
  and $h_{T}(v) = h_{T'}(v)$ for each vertex $v$ of $U_{j}$ with $1 \leq j <
  i$.
  
  The order $<_{lex}$ is a linear extension of $<_{height}$, {\it i. e.} given 
  two tilings $T$ and $T'$ such that $T  <_{height} T'$, we have $T  
  <_{lex} T'$.
  
\begin{proposition}
   Let $T$ be a (non-maximal) tiling of $F$ and let $T_{succ}$ denote the 
   successor of $T$ in the lexicographic order;  
 let $i$ denote the largest integer such that an upward flip is 
 possible in $U_{i}$. 
 
 The tiling  $T_{succ}$ is the lowest tiling (for 
 $<_{height}$) such that 
 $h_{T_{succ}}(v) = h_{T}(v) +4$ for each vertex $v$ of $U_{i}$,
 and $h_{T}(v) = h_{T_{succ}}(v)$ for
  each vertex $v$ of $U_{j}$ with $1 \leq j < i$.

\end{proposition}

\begin{proof}
  Let $T'$ be a tiling such that $h_{T_{succ}}(v) = h_{T}(v) +4$ for each
  vertex $v$ of $U_{i}$ and $h_{T}(v) = h_{T_{succ}}(v)$ for each vertex $v$ 
  of $U_{j}$ with $1 \leq j < i$. 
  By definition, one has $T_{succ} \leq_{lex} T'$.
  
  Moreover, assume that $h_{T} = h_{T_{succ}}$ in $U_{i}$. Thus, by 
  Corollary \ref{c:isomorphic-height-flip}, one can pass from $T$ to 
  $sup(T, T_{succ})$ by a sequence of upward flips, which contradicts 
  the definition of the integer $i$. 
\end{proof}

This proposition enables one to generate all the tilings of $F$ as 
follows:

{\bf Initialization:} Construct the graph $G_{forc.}$ of forced 
components, the tiling 
$T_{min}$, and $o(T_{min})$. Output the tiling $T_{min}$.

The algorithm uses a variable tiling $T$ stored in memory, 
which for initialization is 
equal to $T_{min}$.

{\bf Main loop:} 
Compute the successor of 
$T$ as follows: 

\begin{itemize}
\item  Find the last component $U_{i}$ on which an upward flip can be 
done (if no upward flip is possible, then stop).

\item  Construct the minimal tiling $T$' such that for each vertex $v$ 
of $\cup_{j =1}^{i-1} U_{j}$, $h_{T}(v) = h_{T'}(v)$ and for each vertex $v$ 
of $ U_{i}$, $h_{T'}(v) = h_{T}(v)+4$.

\item  Replace $T$ by $T'$, $o(T)$ by $o(T')$, output the tiling and go 
back to the beginning of the loop. 
\end{itemize}

The second item of the main loop can be done in $O(n^{2})$ time units using an 
algorithm derived from the algorithm of 
construction of the minimal tiling (it suffices to change the 
initialization, fixing appropriate value of $h(v)$ for $v$ in 
$\cup_{j =1}^{i} U_{j}$). 

Thus, once the initialization is 
done, the maximal waiting time between two consecutive tilings is  
$O(n^{2})$ time units. The memory space is $O(n^2)$ since for 
each vertex $v$, one has to store $h_{T}(v)$ (using unary numeration).

\subsection{Uniform random sampling}

Consider the  following process: Given a tiling $T$, find at random a 
forced component  $C$ and a direction (upwards or downwards). If a
flip can be done in $C$ according to the chosen direction, then make 
this flip; otherwise, do not change $T$. 
Trivially, this Markovian random process is ergodic and converges to the 
uniform distribution.

Moreover, the method of ``coupling from the past'' \cite{Prop-Wil} can be
applied since the process is monotonic and one has a method to construct the
maximal and minimal tilings. We thus have a randomized algorithm to sample
domino tilings uniformly at random. The space required is polynomial. 

It has been previously proved \cite{lubyrandallsinclair} \cite{Wil} that 
this process is rapidly mixing for figures without holes. We 
conjecture that is remains true in the general case.

About algorithms, the reader can also easily verify that,  
given a pair $(T, T')$ of tilings, one can 

\begin{itemize} 
\item compute in linear time if 
one can pass from $T$ to $T'$ by a sequence of local flips (it suffices 
to compare $h_{T}$ and $h_{T'}$ on the boundary of $F$);

\item  compute in polynomial time with (low degree)  a shortest path of 
 flips to go from $T$ to $T'$ (using the cyclic orientations and the 
 distance) and the length of such a path.  

\end{itemize}

\section{Extension to other types of  tilings}

\subsection{Calissons}
The same study can easily be done for calisson  ({\it i. e.} tiles 
formed by two neighbouring cells of the triangular lattice) tilings to get 
similar results. In this case, local flips are induced by the two
tilings of hexagons formed by six triangular cells.  There are only 
 two small differences, detailed below:  

\begin{itemize} 
\item  One has two types of connectivity for triangular cells 
(3-connectivity  for cells which share an edge, and 12-connectivity 
for cells which share a vertex). Thus, some vertices have to be 
duplicated or triplicated (see Figure \ref{triplic}). 

\begin{figure}[htbp]
\centerline{\epsfig{file=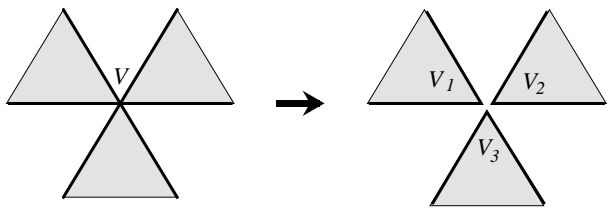,width=8cm}}
\caption{Example of ``triplication''.} 
\label{triplic}
\end{figure}

\item For the proof of corollary \ref{coro cycle critique}, one has to 
consider a part of the critical cycle with vertices $v = (x, y)$ such 
that $y$ is maximal. 

\begin{figure}[htbp]
\centerline{\epsfig{file=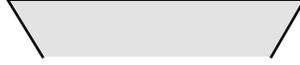,width=4cm}}
\caption{A part of a critical cycle which implies that there is 
no tiling} 
\label{triangle_critique}
\end{figure}

\end{itemize} 

\subsection{Bicolored Wang tiles}

The case of dominoes is a particular case of tilings with Wang 
tiles ({\it i. e.}   $1 \times 1$ squares with colored edges, see 
\cite{Moo-Rap-Rem} for details).  They give rise to a 
tiling if  the colors on the edges of neighbour squares are compatible. 

An instance of the problem of tiling by bicolored Wang tiles is given
by a finite  figure and
a coloration of the edges which are on its boundary. Hence, a domino 
tiling is a tiling by Wang tiles with one red edge and three blue 
edges, all the edges on the boundary of the figure being blue.

\subsubsection{Eulerian orientations}
The same study can  be done for tilings with ``balanced  Wang tiles''
such that each square has two blue edges and two red edges. In this
case, one has $hd_{T}(v, v') = sp(v, v')+eq(v,v')$ if the corresponding edge 
is blue and $hd_{T}(v, v') = -sp(v, v')+eq(v,v')$ otherwise. 
 This is
recognizable as the height function for Eulerian orientations of the
dual lattice, called the {\em six-vertex ice model} by physicists
\cite{beijeren}.  This is also equivalent to the height function for
three-colorings of vertices of the square lattice, and to alternating-sign
matrices \cite{proppdcg}. The results are similar to those obtained for 
dominoes.

\subsubsection{Examples with finite height functions}

 For the case of ``odd tiles'' ({\it i. e.} tiles with exactly three  edges of 
 the same color (blue or red), see 
\cite{Moo-Rap-Rem}): One has to take a height function in $\Z/8\Z$ 
such that  $hd_{T}(v, v') = sp(v, v')+eq(v,v')$ if the corresponding edge 
is blue and $hd_{T}(v, v') = -3sp(v, v') +eq(v,v')$ otherwise. With our
technique of 
equilibrium value, it is easily proved that the set of the tilings of a 
fixed figure has a structure of boolean lattice (or hypercube), even 
if the figure has  holes. 

The case of ``even tiles'' ({\it i. e.} tiles with an even 
number of blue edges and an even 
number of red edges) is very similar, with
$val_{T}(v, v') = (sp(v, v'), 0) $ if the corresponding edge 
is blue and $val_{T}(v, v') = (sp(v, v'),sp(v, v')) $ otherwise.  
These values are taken in $\Z/4\Z \times \Z/2\Z$.

\end{document}